%
%
%
\newif\ifsect\newif\iffinal
\secttrue\finalfalse
\def\thm #1: #2{\medbreak\noindent{\bf #1:}\if(#2\thmp\else\thmn#2\fi}
\def\thmp #1) { (#1)\thmn{}}
\def\thmn#1#2\par{\enspace{\sl #1#2}\par
        \ifdim\lastskip<\medskipamount \removelastskip\penalty 55\medskip\fi}
\def\square{{\msam\char"03}}
\def\qedn{\thinspace\null\nobreak\hfill\square\par\medbreak}
\def\pf{\ifdim\lastskip<\smallskipamount \removelastskip\smallskip\fi
        \noindent{\sl Proof\/}:\enspace}
\def\itm#1{\item{\rm #1}\ignorespaces}

\def\bar#1{\overline{#1}}
%
%
\def\Figuraeps #1 (#2){\message{Figura #1}
	\midinsert      
	\centerline{\BoxedEPSF{#2.eps}}
	\bigskip
	\centerline{\bf Figure~#1}
	\endinsert}
\def\Figurascaledeps #1 (#2)(#3){\message{Figura #1}
	\midinsert      
	\centerline{\BoxedEPSF{#2.eps scaled #3}}
	\bigskip
	\centerline{\bf Figure~#1}
	\endinsert}
\newcount\parano
\newcount\eqnumbo
\newcount\thmno
\newcount\versiono
\newcount\remno
\newbox\notaautore
\def\neweqt#1$${\xdef #1{(\number\parano.\number\eqnumbo)}
    \eqno #1$$
    \global \advance \eqnumbo by 1}
\def\newrem #1\par{\global \advance \remno by 1
    \medbreak
{\bf Remark \the\parano.\the\remno:}\enspace #1\par
\ifdim\lastskip<\medskipamount \removelastskip\penalty 55\medskip\fi}
\def\newthmt#1 #2: #3{\xdef #2{\number\parano.\number\thmno}
    \global \advance \thmno by 1
    \medbreak\noindent
    {\bf #1 #2:}\if(#3\thmp\else\thmn#3\fi}
\def\neweqf#1$${\xdef #1{(\number\eqnumbo)}
    \eqno #1$$
    \global \advance \eqnumbo by 1}
\def\newthmf#1 #2: #3{\xdef #2{\number\thmno}
    \global \advance \thmno by 1
    \medbreak\noindent
    {\bf #1 #2:}\if(#3\thmp\else\thmn#3\fi}
\def\forclose#1{\hfil\llap{$#1$}\hfilneg}
\def\newforclose#1{
	\ifsect\xdef #1{(\number\parano.\number\eqnumbo)}\else
	\xdef #1{(\number\eqnumbo)}\fi
	\hfil\llap{$#1$}\hfilneg
	\global \advance \eqnumbo by 1}
\def\forevery#1#2$${\displaylines{\let\eqno=\forclose
        \let\neweq=\newforclose\hfilneg\rlap{$\qquad\quad\forall#1$}\hfil#2\cr}$$}
\def\noNota #1\par{}
\def\today{\ifcase\month\or
   January\or February\or March\or April\or May\or June\or July\or August\or
   September\or October\or November\or December\fi
   \space\number\year}
\def\inizia{\ifsect\let\neweq=\neweqt\else\let\neweq=\neweqf\fi
\ifsect\let\newthm=\newthmt\else\let\newthm=\newthmf\fi}
\def\bititolo{\empty}
\gdef\begin #1 #2\par{\xdef\titolo{#2}
\ifsect\let\neweq=\neweqt\else\let\neweq=\neweqf\fi
\ifsect\let\newthm=\newthmt\else\let\newthm=\newthmf\fi
\iffinal\let\Nota=\noNota\fi
\centerline{\titlefont\titolo}
\if\bititolo\empty\else\medskip\centerline{\titlefont\bititolo}
\xdef\titolo{\titolo\ \bititolo}\fi
\bigskip
\centerline{\bigfont
\autore \ifvoid\notaautore\else\footnote{${}^1$}{\unhbox\notaautore}\fi}
\bigskip\if\istituto!\centerline{\today}\else
\centerline{\istituto}
\centerline{\indirizzo}
\centerline{E-mail: \email}
\medskip
\centerline{#1\ \anno}\fi
\bigskip\bigskip
\ifsect\else\global\thmno=1\global\eqnumbo=1\fi}
\def\anno{2007}
\font\titlefont=cmssbx10 scaled \magstep1
\font\bigfont=cmr12
\font\eightrm=cmr8
\font\sc=cmcsc10
\font\bbr=msbm10
\font\sbbr=msbm7
\font\ssbbr=msbm5
\font\msam=msam10

\font\bfm=cmmib10

\def\ca #1{{\cal #1}}

\nopagenumbers
\binoppenalty=10000
\relpenalty=10000
\newfam\amsfam
\textfont\amsfam=\bbr \scriptfont\amsfam=\sbbr \scriptscriptfont\amsfam=\ssbbr
\newfam\boldifam
\textfont\boldifam=\bfm
\let\de=\partial
\let\eps=\varepsilon

\def\Re{\mathop{\rm Re}\nolimits}

\mathchardef\void="083F
\mathchardef\ellb="0960
\mathchardef\taub="091C
\def\C{{\mathchoice{\hbox{\bbr C}}{\hbox{\bbr C}}{\hbox{\sbbr C}}
{\hbox{\sbbr C}}}}
\def\R{{\mathchoice{\hbox{\bbr R}}{\hbox{\bbr R}}{\hbox{\sbbr R}}
{\hbox{\sbbr R}}}}
\def\N{{\mathchoice{\hbox{\bbr N}}{\hbox{\bbr N}}{\hbox{\sbbr N}}
{\hbox{\sbbr N}}}}

\newcount\notitle
\notitle=1
\headline={\ifodd\pageno\rhead\else\lhead\fi}
\def\rhead{\ifnum\pageno=\notitle\iffinal\hfill\else\hfill\tt Version
\the\versiono; \the\day/\the\month/\the\year\fi\else\hfill\eightrm\titolo\hfill
\folio\fi}
\def\lhead{\ifnum\pageno=\notitle\hfill\else\eightrm\folio\hfill\autore\hfill
\fi}
\newbox\bibliobox
\def\setref #1{\setbox\bibliobox=\hbox{[#1]\enspace}
    \parindent=\wd\bibliobox}
\def\biblap#1{\noindent\hang\rlap{[#1]\enspace}\indent\ignorespaces}
\def\art#1 #2: #3! #4! #5 #6 #7-#8 \par{\biblap{#1}#2: {\sl #3\/}.
    #4 {\bf #5} (#6)\if.#7\else, \hbox{#7--#8}\fi.\par\smallskip}
\def\book#1 #2: #3! #4 \par{\biblap{#1}#2: {\bf #3.} #4.\par\smallskip}
\def\coll#1 #2: #3! #4! #5 \par{\biblap{#1}#2: {\sl #3\/}. In {\bf #4,}
#5.\par\smallskip}
\def\pre#1 #2: #3! #4! #5 \par{\biblap{#1}#2: {\sl #3\/}. #4, #5.\par\smallskip}
%

%

\def\Nota #1\par{}
\newcount\defno
\def\smallsect #1. #2\par{\bigbreak\noindent{\bf #1.}\enspace{\bf #2}\par
    \global\parano=#1\global\eqnumbo=1\global\thmno=1\global\defno=0\global\remno=0
    \nobreak\smallskip\nobreak\noindent\message{#2}}
\def\newdef #1\par{\global \advance \defno by 1
    \medbreak
{\bf Definition \the\parano.\the\defno:}\enspace #1\par
\ifdim\lastskip<\medskipamount \removelastskip\penalty 55\medskip\fi}
\def\autore{Marco Abate\footnote{${}^1$}{\sevenrm Dipartimento di Matematica, Universit\`a di Pisa, Largo Pontecorvo 5, 56127 Pisa, Italy, e-mail: abate@dm.unipi.it}, 
Alberto Saracco\footnote{${}^2$}{\sevenrm Dipartimento di Matematica, Universit\`a di Parma, Viale G.P. Usberti 53/A, 43124 Parma, Italy, e-mail: alberto.saracco@unipr.it}}
\def\istituto{!}
\def\indirizzo{${}^2$ }
\def\email{}
\def\bititolo{in strongly pseudoconvex domains}
\year=2009
\finaltrue
\versiono=11
\begin {September} Carleson measures and uniformly discrete sequences

{\narrower
{\sc Abstract.} 
We characterize using the Bergman kernel Carleson measures of Bergman spaces in strongly pseudoconvex bounded
domains in~$\C^n$, generalizing to this setting theorems proved by Duren and Weir
for the unit ball.  We also show that uniformly discrete (with respect to the Kobayashi
distance) sequences give examples of Carleson measures, and we compute the speed
of escape to the boundary of uniformly discrete sequences in strongly pseudoconvex
domains, generalizing results obtained in the unit ball by Jevti\'c, Massaneda and Thomas, by Duren and Weir, and by MacCluer.

}

\smallsect 0. Introduction

In his celebrated solution of the corona problem in the disk, Carleson [C] introduced an important class of measures to study the structure of the Hardy spaces
\def\autore{Marco Abate, Alberto Saracco}
of the unit disk $\Delta\subset\C$. Let $A$ be a Banach space of holomorphic functions 
on a domain~$D\subset\C^n$, and assume that $A$ is contained in $L^p(D)$ for some~$p>0$. A finite positive Borel measure $\mu$ 
on~$D$ is a {\sl Carleson measure} of $A$ if there exists a constant $C>0$ such that
$$
\forevery{f\in A}\int_D |f|^p\,d\mu\le C\|f\|^p_A\;.
$$
Carleson studied Carleson measures of the Hardy spaces~$H^p(\Delta)$, showing that
a finite positive Borel measure $\mu$ is a Carleson measure of $H^p(\Delta)$ if and only if 
there exists a constant $C>0$ such that
$\mu(S_{\theta_0,h})\le C h$ for all sets
$$
S_{\theta_0,h}=\{re^{i\theta}\in\Delta\mid 1-h\le r<1,\ |\theta-\theta_0|\le h\}
$$
(see [D]); in particular the set of Carleson measures of $H^p(\Delta)$ does not depend
on~$p$.

In 1975, Hastings [H] (see also Oleinik and Pavlov [OP] and Oleinik [O]) proved a similar characterization for the Carleson measures of the Bergman spaces~$A^p(\Delta)$: a finite positive Borel measure $\mu$ is a Carleson measure of $A^p(\Delta)$
if and only if there exists a constant $C>0$ such that $\mu(S_{\theta_0,h})\le Ch^2$  for all
$\theta\in[0,2\pi]$ and $h\in(0,1)$. As a consequence, again, the set of Carleson measures of $A^p(\Delta)$ does not depend on~$p$. 

The sets $S_{\theta_0,h}$ clearly are not invariant under automorphisms of the disk, whereas
one would like to characterize Carleson measures for Bergman spaces in an invariant way,
related to the intrinsic (hyperbolic) geometry of the disk and not to the extrinsic Euclidean geometry. 
This has been done in 1983 by Luecking [Lu]; indeed he proved (see also [DS, Theorem~14, p.~62]) that a finite positive Borel measure~$\mu$ is a Carleson measure of~$A^p(\Delta)$ 
if and only if for some (and hence all) $0<r<1$ there is a constant $C_r>0$ such that
$\mu\bigl(B_\Delta(z_0,r)\bigr)\le C_r\nu\bigl(B_\Delta(z_0,r)\bigr)$ for all $z_0\in\Delta$, where 
$\nu$ is the  Lebesgue (area) measure and $B_\Delta(z_0,r)\subset\Delta$ is the Poincar\'e disk
$$
B_\Delta(z_0,r)=\left\{z\in\Delta\biggm| \left|{z-z_0\over 1-\bar{z_0}z}\right|<r\right\}
$$
of center~$z_0\in\Delta$ and pseudohyperbolic radius~$r$. 

The first characterization for Carleson measures of the Bergman spaces of the unit ball $B^n\subset\C^n$ has been
given by Cima and Wogen [CW], using again sets defined in terms of the Euclidean geometry
of~$B^n$, and thus not invariant under automorphisms. Again, as essentially noticed by
Luecking [Lu] and explicitely stated by Duren and Weir [DW], it is possible to give 
a characterization for the Carleson measures of Bergman spaces of~$B^n$ by using 
the balls for the Bergman (or Kobayashi, or pseudohyperbolic) distance: a  
finite positive Borel measure~$\mu$ is a Carleson measure of~$A^p(B^n)$ 
if and only if for some (and hence all) $0<r<1$ there is a constant $C_r>0$ such that
$\mu\bigl(B_{B^n}(z_0,r)\bigr)\le C_r\nu\bigl(B_{B^n}(z_0,r)\bigr)$ for all $z_0\in B^n$, where 
$\nu$ is the  Lebesgue $2n$-dimensional measure and $B_{B^n}(z_0,r)\subset B^n$ is the ball
for the Bergman distance of center $z_0\in B^n$ and radius~${1\over 2}\log{1+r\over 1-r}\in\R^+$
(that is, of radius $r$ in the pseudohyperbolic distance; recall that the pseudohyperbolic
distance~$\rho$ is related to the Bergman or Kobayashi distance~$k_{B^n}$ by 
the formula $\rho=\tanh(k_{B^n})$).

In 1995, Cima and Mercer [CM] characterized Carleson measures of the Bergman spaces
of strongly pseudoconvex domains; in particular, they proved that in this case too the class
of Carleson measures of~$A^p(D)$ does not depend on~$p$. Their characterization is stated again in terms of the extrinsic Euclidean geometry of the domain, but the proof uses in an essential way the intrinsic geometry of strongly convex domains, as well as the
construction of particular {\it ad hoc} functions. 

A particularly important Bergman space is, of course, $A^2(D)$, where the Bergman
kernel lives. This suggests the question of whether it is possible to characterize Carleson
measures using the Bergman kernel. This has been done by Duren and Weir [DW]
for the unit ball; our first main result is the generalization of their characterization 
to strongly pseudoconvex domains.

Let $K\colon D\times D\to\C$ be the Bergman kernel of a strongly pseudoconvex domain
$D\subset\subset\C^n$. For any finite positive Borel measure $\mu$ on~$D$, the {\sl Berezin
transform} of~$\mu$ is the function $B\mu\colon D\to\R$ given by
$$
B\mu(z)=\int_D {|K(\zeta,z)|^2\over K(z,z)}\,d\mu(\zeta)\;.
$$
Then (see Theorem~2.4 for a more complete statement):

\newthm Theorem \iCarleson: Let $\mu$ be a finite positive Borel measure on a strongly pseudoconvex bounded  domain $D\subset\subset\C^n$. Then $\mu$ is a Carleson measure of $A^p(D)$ if and only if its Berezin transform~$B\mu$ is bounded.

Previous proofs of characterizations of Carleson measures in the ball 
heavily relied on the homogeneity of the ball under its automorphisms
group, and on the explicit expression of the automorphisms; but this approach cannot
be used in our setting, because strongly pseudoconvex domains not biholomorphic to the ball admit very few automorphisms. Our proof depends instead on
a detailed understanding of the 
intrinsic Kobayashi geometry of strongly pseudoconvex domains, and on Fefferman's 
estimates on the Bergman kernel.

A natural question is how to construct explicit examples of Carleson measures
in strongly pseudoconvex 
domains. As in the unit disk and in the unit ball, an important family of examples is
provided by uniformly discrete sequences. Let $(X,d)$ be a metric space; a sequence 
$\{x_j\}\subset X$ of points in $X$ is {\sl uniformly discrete} if there exists $\delta>0$
such that $d(x_j,x_k)\ge\delta$ for all $j\ne k$. Then we are able to prove the
following generalization of results proved by Massaneda [Ma], Jevti\'c, Massaneda and
Thomas [JMT] and Duren and Weir [DW] in the unit ball (see Theorem~3.2 for a more
complete statement):

\newthm Theorem \zdtre: Let $D\subset\subset\C^n$ be a strongly pseudoconvex bounded  domain, and let $\Gamma=\{z_j\}$ be a sequence
in~$D$. Then $\Gamma$ is a finite union of uniformly discrete (with respect to the Kobayashi
distance) sequences if and only if 
$ \sum\limits_{z_j\in\Gamma} d(z_j,\de D)^{n+1}\delta_{z_j}$
is a Carleson measure of $A^p(D)$, where $\delta_{z_j}$ is the Dirac measure at~$z_j$
and $d(\cdot,\de D)$ is the Euclidean distance from the boundary.

Finally, Duren, Schuster and Vukoti\'c [DSV], Duren and Weir [DW], and MacCluer [Mc] also studied
how fast in the unit ball a uniformly discrete (with respect to
the pseudohyperbolic or Bergman distances of the ball) sequence
escapes to the boundary. We are able to
generalize their results to strongly pseudoconvex domains (see Theorem~3.6, and
Proposition~3.4 for a similar result not requiring pseudoconvexity): 

\newthm Theorem \zdqua: Let $D\subset\C^n$ be a strongly pseudoconvex bounded  domain. 
Let $\Gamma=\{z_j\}\subset D$ be an uniformly discrete (with respect to the Kobayashi distance) sequence with $d(z_j,\de D)<1$ for all~$j\in\N^*$. Then
$$
\sum_{z_j\in\Gamma} d(z_j,\de D)^{n}\, h\left(-{1\over\log d(z_j,\de D)}\right)<+\infty
\neweq\eqzuno
$$
for all increasing functions $h\colon\R^+\to\R^+$ such that
$$
\sum_{m=1}^{+\infty}h\left({1\over m}\right)<+\infty\;.
$$


\smallsect 1. The intrinsic geometry of strongly convex domains 

In this section we shall prove a number of estimates on the intrinsic geometry
of strongly convex domains, as described by the Kobayashi distance.
In particular, we shall study the boundary behavior of Kobayashi balls,
and we shall prove a sort of submean property for nonnegative plurisubharmonic
functions in Kobayashi balls.

Let us briefly recall the definition and the main properties of the Kobayashi distance;
we refer to [A], [JP] and [K] for details and much more. Let $k_\Delta$ denote the 
Poincar\'e distance on the unit disk~$\Delta\subset\C$. If $X$ is a complex manifold,
the {\sl Lempert function}~$\delta_X\colon X\times X\to\R^+$ of~$X$ is defined by
$$
\delta_X(z,w)=\inf\{k_\Delta(\zeta,\eta)\mid
\hbox{there exists a holomorphic $\phi\colon\Delta\to X$ with $\phi(\zeta)=z$ and
$\phi(\eta)=w$}\}
$$
for all $z$, $w\in X$. The {\sl Kobayashi pseudodistance}~$k_X\colon X\times X\to\R^+$ of~$X$
is the smallest pseudodistance on~$X$ bounded below by~$\delta_X$. We say that $X$ is
{\sl (Kobayashi) hyperbolic} if $k_X$ is a true distance --- and in that case it is known that the
metric topology induced by~$k_X$ coincides with the manifold topology of~$X$ (see, e.g.,
[A, Proposition~2.3.10]). For instance, all bounded domains are hyperbolic (see, e.g., [A, Theorem~2.3.14]).

The main property of the Kobayashi (pseudo)distance is that it is contracted by holomorphic maps: if $f\colon X\to Y$ is a holomorphic map then
$$
\forevery{z,w\in X} k_Y\bigl(f(z),f(w)\bigr)\le k_X(z,w)\;.
$$
In particular, the Kobayashi distance is invariant under biholomorphisms, and decreases under inclusions:
if $D_1\subset D_2\subset\subset\C^n$ are two bounded domains we have 
$k_{D_2}(z,w)\le k_{D_1}(z,w)$ for all $z$,~$w\in D_1$.

It is easy to see that the Kobayashi distance of the unit disk coincides with the Poincar\'e distance.
Furthermore, the Kobayashi distance of the unit ball $B^n\subset\C^n$ coincides with the
Bergman distance (see, e.g., [A, Corollary~2.3.6]). 

If $X$ is a hyperbolic manifold, $z_0\in X$ and $r\in(0,1)$ we shall denote by~$B_X(z_0,r)$
the Kobayashi ball of center~$z_0$ and radius ${1\over2}\log{1+r\over1-r}$:
$$
B_X(z_0,r)=\{z\in X\mid \tanh k_X(z_0,z)<r\}\;.
$$
Notice that $\rho_X=\tanh k_X$ is still a distance on~$X$, because $\tanh$ is a strictly convex function on~$\R^+$. In particular, $\rho_{B^n}$ is the pseudohyperbolic distance of~$B^n$.

The Kobayashi distance of bounded strongly pseudoconvex domains enjoys several important 
properties. First of all, it is complete (see, e.g., [A, Corollary~2.3.53]), and hence
closed Kobayashi balls are compact. Furthermore, we can very precisely
describe the boundary behavior of the Kobayashi distance: if $D\subset\subset\C^n$
is a strongly pseudoconvex bounded domain and $z_0\in D$, there exist
$c_0$, $C_0>0$ such that
$$
\forevery{z\in D} c_0-{\textstyle{1\over2}}\log d(z,\de D)\le k_D(z_0,z)\le 
C_0-{\textstyle{1\over2}}\log d(z,\de D)\;,
\neweq\equbb
$$
where $d(\cdot,\de D)$ denotes the Euclidean distance from the boundary of~$D$
(see [A, Theorems~2.3.51 and~2.3.52]).

Let us finally recall a couple of facts on Kobayashi balls of~$B^n$; for proofs see [A, Section~2.2.2], [Ru, Section~2.2.7] and [DW]. The ball $B_{B^n}(z_0,r)$ is given by
$$
B_{B^n}(z_0,r)=\left\{z\in B^n\biggm|{(1-\|z_0\|^2)(1-\|z\|^2)\over|1-\langle z,z_0\rangle|^2}> 1-r^2\right\}\;.
\neweq\eqeqball
$$
Geometrically, it is an ellipsoid
of (Euclidean) center
$$
c={1-r^2\over 1-r^2\|z_0\|^2}\,z_0\;,
$$
its intersection with the complex line $\C z_0$ is an Euclidean disk
of radius 
$$
r\,{1-\|z_0\|^2\over 1-r^2\|z_0\|^2}\;,
$$
and its intersection with the affine subspace through~$z_0$
orthogonal to~$z_0$ is an Euclidean ball of the larger radius 
$$
r\,\sqrt{{1-\|z_0\|^2\over1-r^2\|z_0\|^2}}\;.
$$

Let $\nu$ denote the Lebesgue volume measure of~$\R^{2n}$,
normalized so that $\nu(B^n)=1$. Then the volume of a Kobayashi ball $B_{B^n}(z_0,r)$ is given by
(see [DW])
$$
\nu\bigl(B_{B^n}(z_0,r)\bigr)=
r^{2n}\left({1-\|z_0\|^2\over 1-r^2\|z_0\|^2}\right)^{n+1}\;.
\neweq\eqvolkba
$$

A similar estimate holds for the volume of Kobayashi balls in
strongly pseudoconvex bounded domains:

\newthm Lemma \sei: Let $D\subset\subset\C^n$ be a strongly pseudoconvex bounded  domain. Then there exist $c_1>0$ and, for each $r\in(0,1)$, a $C_{1,r}>0$ depending on~$r$ such that
$$
c_1 r^{2n} d(z_0,\de D)^{n+1}\le \nu\bigl(B_D(z_0,r)\bigr)\le
C_{1,r}
d(z_0,\de D)^{n+1}
$$
for every $z_0\in D$ and $r\in(0,1)$.

\pf 
%
%
Since $D$ is bounded and with smooth boundary, 
there is $\delta>0$ so that any euclidean ball internally tangent to $\de D$ with radius~$\delta$ is completely contained in~$D$. Take $z_0\in D$ with $d(z_0,\de D)<\delta$, and choose $x\in\de D$ so that $d(z_0,\de D)=\|x-z_0\|$. Then the euclidean ball~$B$
tangent to~$\de D$ in~$x$ and with radius~$\delta$ is contained in~$D$, contains~$z_0$
and $d(z_0,\de D)=d(z_0,\de B)$. Since the Kobayashi distance decreases under inclusions, $B_B(z_0,r)\subset B_D(z_0,r)$.
Hence \eqvolkba\ yields (assuming without loss of generality that $B$ is centered at the origin)
$$
\eqalign{
\nu\bigl(B_D(z_0,r)\bigr)&\ge\nu\bigl(B_B(z_0,r)\bigr)=\delta^{2n}\nu\bigl(B_{B^n}(z_0/\delta,r)
\bigr)\cr
&=
\delta^{2n}r^{2n}\left({1-\|z_0/\delta\|^2\over1-r^2\|z_0/\delta\|^2} \right)^{n+1}
=
\delta^{2n}r^{2n}\left({\delta^2-\|z_0\|^2\over\delta^2-r^2\|z_0\|^2} \right)^{n+1}\cr
&\ge
\delta^{2n}r^{2n}\left({\delta(\delta-\|z_0\|)\over\delta^2}\right)^{n+1}
=
\delta^{n-1}r^{2n}d(z_0,\de D)^{n+1}\;,
\cr}
$$
and we are done in this case.

\Nota Another estimate, if needed, is
$$
\nu\bigl(B(z_0,r)\bigr)\ge\nu(B^n)\delta^{2n}r^{2n}\left(1+{\delta(1-r^2)\over d(z_0,\de D)}
\right)^{-(n+1)}\;;
$$
it follows writing $\delta^2-r^2\|z_0\|^2=\delta^2-\|z_0\|^2+\|z_0\|^2(1-r^2)$.

Assume now $d(z_0,\de D)\ge\delta$, and let $B\subset D$ the euclidean ball
of center~$z_0$ and radius~$d(z_0,\de D)$. Then $B_B(z_0,r)$ is an euclidean ball 
of radius $r d(z_0,\de D)$, and so
$$
\nu\bigl(B_D(z_0,r)\bigr)\ge\nu\bigl(B_B(z_0,r)\bigr)=
r^{2n}d(z_0,\de D)^{2n}
\ge 
\delta^{n-1}r^{2n} d(z_0,\de D)^{n+1}\;,
$$ 
and we have obtained the lower estimate in this case too.

For the upper estimate, results of [KM] and [Li] show that there is a $\delta_1>0$ such that
if $d(z_0,\de D)<\delta_1$ then there exist $a_1(r)$, $a_2(r)>0$ so that
$B_D(z_0,r)$ is contained in a (possibly rotated) polydisk~$P$ of center~$z_0$
and polyradius $\bigl(a_1(r)d(z_0,\de D), a_2(r)\sqrt{d(z_0,\de D)},\ldots,
a_2(r)\sqrt{d(z_0,\de D)}\bigr)$. Thus if $d(z_0,\de D)<\delta_1$ we have
$$
\nu\bigl(B_D(z_0,r)\bigr)\le\nu(P)\le ca_1(r)^2a_2(r)^{2(n-1)}d(z_0,\de D)^{n+1}
$$
for a suitable constant $c>0$ independent of~$r$. 

Finally, if $d(z_0,\de D)\ge\delta_1$ we have
$$
\nu\bigl(B_D(z_0,r)\bigr)\le{\nu(D)\over\delta_1^{n+1}}d(z_0,\de D)^{n+1}\;,
$$
and we are done.\qedn


\newrem If $D$ is strongly convex, using Euclidean balls externally tangent to~$D$ it is
possible to show that one can take
$$
C_{1,r}=C_1 r^{2n}(1-r^2)^{-(n+1)}\;,
$$
where $C_1>0$ is a constant independent of~$r$.

%
%


The next lemma contains an estimate on the Euclidean size of Kobayashi balls.

\newthm Lemma \cinque: Let $D\subset\subset\C^n$ be a strongly pseudoconvex 
bounded  domain. Then there is $C_2>0$ such that for every $z_0\in D$
and $r\in(0,1)$ one has
$$
\forevery{z\in B_D(z_0,r)}
{C_2\over 1-r}\, d(z_0,\de D)\ge d(z,\de D)\ge {1-r\over C_2}d(z_0,\de D)\;.
$$

\pf Let us fix $w_0\in D$. Then \equbb\ yields $c_0$, $C_0>0$ such that
$$
c_0-{\textstyle{1\over2}}\log d(z,\de D)\le k_D(w_0,z)\le k_D(z_0,z)+k_D(z_0,w_0)\le
{\textstyle{1\over 2}}\log{1+r\over1-r}+C_0-{\textstyle{1\over2}}\log d(z_0,\de D)\;,
$$
for all $z\in B_D(z_0,r)$, and hence
$$
e^{2(c_0-C_0)}d(z_0,\de D)\le {2\over 1-r}d(z,\de D)\;.
$$
The left-hand inequality is obtained in the same way reversing the roles of~$z_0$ and~$z$.\qedn



\newrem Using again externally tangent Euclidean balls it is possible to show that if $D$ is strongly convex then one can take $C_{2}=4$.

The two previous lemmas (together with the following Corollary~1.7) give the main properties needed
in Luecking's approach [Lu] to the characterization of Carleson measures. However, to get the characterization involving the Berezin transform, we shall need precise information on the
behaviour of the Bergman kernel in Kobayashi balls, that we shall deduce from
Fefferman's estimates using another estimate on the shape of Kobayashi balls. 
In the unit ball the latter estimate has the following form:

\newthm Lemma \bisuno: Let $B^n\subset\C^n$ be the unit ball in~$\C^n$, and
take $z_0\in B^n$ and $r\in(0,1)$. Then
$$
\forevery{z\in B_{B^n}(z_0,r)}1-\|z_0\|^2>{1-r^2\over 4}\left(
\|z-z_0\|^2+|\langle z-z_0,z_0\rangle|\right)\;.
$$

\pf First of all, it is easy to check that $z\in B_{B^n}(z_0,r)$ if and
only if
$$
{(1-\|z_0\|^2)(1-\|z\|^2)\over|1-\langle z,z_0\rangle|^2}>1-r^2\;;
$$
therefore it suffices to prove that
$$
{|1-\langle z,z_0\rangle|^2\over 1-\|z\|^2}\ge{1\over 4}
\left(\|z-z_0\|^2+|\langle z-z_0,z_0\rangle|\right)\;.
\neweq\eqbisuno
$$

Let us write $z_0=\lambda z+w$, with $\lambda\in\C$ and $w\perp z$. Then
$\langle z,z_0\rangle=\bar\lambda\|z\|^2$ and $\|z-z_0\|^2=|1-\lambda|^2\|z\|^2
+\|w\|^2$; moreover $|\lambda|^2\|z\|^2+\|w\|^2=\|z_0\|^2<1$, that is
$\|w\|^2<1-|\lambda|^2\|z\|^2$. Now, a computation shows that
$\|z\|^4|\lambda-1|^2\ge 0$ is equivalent to
$$
|1-\bar\lambda\|z\|^2|^2\ge(1-\|z\|^2)\bigl(|1-\lambda|^2\|z\|^2+1-|\lambda|^2
\|z\|^2\bigr)\;;
$$ 
therefore
$$
{|1-\langle z,z_0\rangle|^2\over 1-\|z\|^2}={|1-\bar\lambda\|z\|^2|^2\over
1-\|z\|^2}\ge |1-\lambda|^2\|z\|^2+1-|\lambda|^2\|z\|^2>
|1-\lambda|^2\|z\|^2+\|w\|^2=\|z-z_0\|^2\;.
\neweq\eqbisdue
$$

Write $z=\mu z_0+w$, with $\mu\in\C$ and $w\perp z_0$. Then
$\langle z,z_0\rangle=\mu\|z_0\|^2$ and $|\langle z-z_0,z_0\rangle|=|1-\mu|\|z_0\|^2$; moreover $|\mu|^2\|z_0\|^2+\|w\|^2=\|z\|^2<1$, implying
$|\mu|\|z_0\|<1$. Now,
$$
{|1-\langle z,z_0\rangle|^2\over1-\|z\|^2}={|1-\mu\|z_0\|^2|^2\over
1-|\mu|^2\|z_0\|^2-\|w\|^2}\ge{|1-\mu\|z_0\|^2|^2\over
1-|\mu|^2\|z_0\|^2}={|1-\mu\|z_0\|^2|\over
1-|\mu|\|z_0\|}\cdot{|1-\mu\|z_0\|^2|\over
1+|\mu|\|z_0\|}\;.
\neweq\eqbistre
$$
By assumption we have $2\ge|\mu|\|z_0\|(1+\|z_0\|)$. A computation shows 
that this is equivalent to 
$$
|1-|\mu|\|z_0\|^2|^2\ge (1-|\mu|\|z_0\|)^2\;;
$$
therefore
$$
{|1-\mu\|z_0\|^2|\over 1-|\mu|\|z_0\|}\ge{|1-|\mu|\|z_0\|^2|\over 1-|\mu|\|z_0\|}\ge 1\;.
\neweq\eqbisqua
$$
Next, $1+\|z_0\|^2\ge 2\|z_0\|\ge2|\mu|\|z_0\|^2\ge 2\|z_0\|^2\Re\mu$; multiplying this by $1-\|z_0\|^2$ we end up with
$$
|1-\mu\|z_0\|^2|^2\ge|1-\mu|^2\|z_0\|^4=|\langle z-z_0,z_0\rangle|^2\;.
\neweq\eqbiscin
$$ 
Since $1+|\mu|\|z_0\|<2$, putting together \eqbistre, \eqbisqua\ and \eqbiscin\ we get
$$
{|1-\langle z,z_0\rangle|^2\over1-\|z\|^2}>{1\over 2}|\langle z-z_0,z_0\rangle|\;.
\neweq\eqbissei
$$
Putting together \eqbisdue\ and \eqbissei\ we get \eqbisuno, and thus
the assertion.\qedn


A {\sl defining function} for a smooth domain
$D\subset\subset\C^n$ is a smooth function $\psi\colon\C^n\to\R$ 
such that $D=\{\psi>0\}$ and the gradient
$\nabla\psi$ is never zero on~$\de D$. If $D$ is strongly pseudoconvex then we can
find a defining function strictly 
plurisubharmonic in a neighbourhood of~$\de D$.

\newthm Lemma \quattrobis: Let $D\subset\subset\C^n$ be a strongly
pseudoconvex bounded  domain, and $\psi\colon\C^n\to\R$ a defining function for~$D$. Then for every $r\in(0,1)$ there exists $c_{2,r}>0$ depending on~$r$ such that for every $z_0\in D$ one has 
$$
\forevery{z\in B_D(z_0,r)}\qquad d(z_0,\de D)\ge c_{2,r}\bigl(
\|z-z_0\|^2+\bigl|\de\psi_{z_0}(z-z_0)|\bigr)
\;.
$$

\pf Fix $r\in(0,1)$. Since $D$ is bounded, it suffices to prove the statement for $z_0$ close enough to~$\de D$.

By Narasimhan's lemma [Kr, Lemma~3.2.3]
we can cover $\de D$ with open sets $U_1,\ldots,U_l\subset\C^n$
so that for each $j=1,\ldots,l$ there is a biholomorphism $\Phi_j\colon U_j\to\Phi_j(U_j)\subset\C^n$ such that $\Phi_j(U_j\cap D)$ is strongly convex. Furthermore, we can assume that 
each $\Phi_j$ is defined in a slightly larger neighborhood, and hence find $c>0$ such 
that $d(z_0,\de D)\ge c d\bigl(\Phi_j(z_0),\de\Phi_j(U_j\cap D)\bigr)$ for all $z_0\in U_j\cap D$
close enough to~$\de D$
and all $j=1,\ldots,l$. For the same reason, and recalling Lemma~\cinque, $|\de\psi_{z_0}(z-z_0)|$ and $\bigl|\de(\psi\circ\Phi_j^{-1})_{\Phi_j(z_0)}\bigl(\Phi_j(z)-\Phi_j(z_0)\bigr)\bigr|$ are comparable as soon as
$z_0\in U_j\cap D$ is close 
enough to~$\de D$.
Finally, by the localization property of the Kobayashi distance
(see [A, Theorem~2.3.65]) if $z_0\in U_j\cap D$ is close enough to the boundary of~$D$
we can find $r_1\in(0,1)$ depending only on~$r$ such that $B_D(z_0,r)\subseteq B_{U_j\cap D}
(z_0,r_1)$. The upshot of these remarks is that it suffices to prove the statement
when $D$ is strongly convex. 

So, assume $D$ strongly convex, and 
let $\delta>0$ be such that if $d(z_0,\de D)<\delta$ then
there exists a unique $x=x(z_0)
\in\de D$ so that $d(z_0,\de D)=\|x-z_0\|$; again it suffices to prove the statement
for $d(z_0,\de D)<\delta$. 

Take then $z_0\in D$ with $d(z_0,\de D)<\delta$; since $D$ is strongly
convex, there exists an euclidean ball $B$ containing~$D$ and tangent
to $D$ in~$x=x(z_0)$; in particular, $d(z_0,\de D)=d(z_0,\de B)$. 
Let $R=R(z_0)>0$ be the radius of~$B$; up to a translation, we can 
assume that $B$ is centered at the origin.

Let $z\in B_D(z_0,r)$. Since $B_D(z_0,r)\subseteq B_B(z_0,r)$, 
Lemma~\bisuno\ implies
$$
2Rd(z_0,\de D)\ge
R^2-\|z_0\|^2>{1-r^2\over 4}\bigl(\|z-z_0\|^2+|\langle z-z_0,z_0\rangle|
\bigr)\;.
$$
Now, $|\langle z-z_0,z_0\rangle|$ is (a uniform multiple of) the distance of~$z$ from the
complex hyperplane~$\pi$ passing through~$z_0$ and parallel to the complex hyperplane tangent to~$\de B$ at~$x$. But the latter coincides with
the complex hyperplane tangent to~$\de D$ at~$x$, which is exactly given
by~$\de\psi_x(z-x)=0$. Therefore $\pi$ has equation $\de\psi_x(z-z_0)=0$. It follows that the difference between $|\langle z-z_0,z_0\rangle|$ and
$|\de\psi_{z_0}(z-z_0)|$ is (uniformly in $z$ and $z_0$) of the order of 
$\|z-z_0\|\|z_0-x\|=d(z_0,\de D)\|z-z_0\|$; so we get
$$
\bigl(2R+C\|z-z_0\|(1-r^2)\bigr)d(z_0,\de D)\ge {1-r^2\over 4}
\bigl(
\|z-z_0\|^2+\bigl|\de\psi_{z_0}(z-z_0)|\bigr)
$$
for a suitable constant $C>0$. Since $D$ is bounded, $R$ too is bounded 
as a function of~$z_0$, and the assertion follows.\qedn

\newrem If $D$ is strongly convex the proof shows that we can take $c_{2,r}=c_2(1-r^2)$
for a suitable $c_2>0$ independent of~$r$.

We now prove a covering lemma for $D$.

\newthm Lemma \uno: Let $D\subset\subset\C^n$ be a strongly pseudoconvex bounded  domain. Then for every $r\in(0,1)$ there exist $m\in\N$ and a sequence 
$\{z_k\}\subset D$ of points such that
$D=\bigcup_{k=0}^\infty B(z_k,r)$
and no point of $D$ belongs to more than $m$ of the balls $B_D(z_k,R)$,
where $R={1\over 2}(1+r)$.

\pf Let $\{B_j\}_{j\in\N}$ be a sequence of Kobayashi balls of radius $r/3$ covering~$D$. We can extract a subsequence $\{\Delta_k=B_D(z_k,r/3)\}_{k\in \N}$ of disjoint balls in the following way: set $\Delta_1=B_1$. Suppose we have already chosen $\Delta_1,\ldots,\Delta_l$. We define $\Delta_{l+1}$ as the first ball in the sequence $\{B_j\}$ which is disjoint from $\Delta_1\cup\cdots\cup\Delta_l$. In particular, by construction every $B_j$
must intersect at least one~$\Delta_k$. 

We now claim that $\{B_D(z_k,r)\}_{k\in\N}$ is a covering of $D$. Indeed, let $z\in D$. Since 
$\{B_j\}_{j\in\N}$ is a covering of~$D$, there is $j_0\in\N$ so that $z\in B_{j_0}$. As remarked above, we get $k_0\in\N$ so that $B_{j_0}\cap\Delta_{k_0}\neq\void$. 
Take $w\in B_{j_0}\cap\Delta_{k_0}$. Then
$$
\rho_D(z,z_{k_0})\leq \rho_D(z,w)+\rho_D(w,z_{k_0})\le {\textstyle{2\over3}}r\;,
$$
and $z\in B_D(z_{k_0},r)$.

To conclude the proof we have to show that there is $m=m_r\in\N$ so that
each point $z\in D$ belongs to at most $m$ of the balls $B(z_k,R)$. Put $r_1={1\over 3}
\min\{r,1-r\}$ and $R_1={1\over 6}(5+r)$. Since $z\in B_D(z_k,R)$ is equivalent to $z_k\in B_D(z,R)$, we have that
$z\in B_D(z_k,R)$ implies $B_D(z_k,r_1)\subset B_D(z,R_1)$. Furthermore, Lemmas~\sei\ 
and~\cinque\ yields
$$
\nu\bigl(B_D(z_k,r_1)\bigr)\ge c_1r_1^{2n}d(z_k,\de D)^{n+1}\ge {c_1\over C_2^{n+1}}
(1-R)^{n+1}r_1^{2n} d(z,\de D)^{n+1}
$$
when $z_k\in B_D(z,R)$. Therefore, since the balls
$B_D(z_k,r_1)$ are pairwise disjoint, using again Lemma~\sei\ we get
$$
\hbox{\rm card}\{k\in\N\mid z\in B_D(z_k,R)\}\le{\nu\bigl(B_D(z,R_1)
\bigr)\over\nu\bigl(B_D(z_k,r_1)\bigr)}\le {C_2^{n+1}C_{1,R_1}\over c_1}{1\over r_1^{2n}(1-R)^{n+1}}\;,
$$
and we are done.\qedn


Our last aim for this section is a sort of submean property in Kobayashi balls for
nonnegative plurisubharmonic functions. Let us first prove it in an
Euclidean ball:

\newthm Lemma \ellisse: Let $B\subset\subset\C^N$ be an Euclidean ball
of radius~$R>0$. Then 
$$
\forevery{z_0\in B\;\forall r\in(0,1)}\qquad
\chi(z_0)\le {4^{n+1}\over R^{n-1}
}{1\over r^{2n} d(z_0,\de B)^{n+1}}\int_{B_B(z_0,r)}\chi\,d\nu
$$
for all nonnegative plurisubharmonic functions $\chi\colon B\to\R^+$. 

\pf Without loss of generality we can assume that $B$ is centered at the
origin. Let $\phi_{z_0/R}\colon B^n\to B^n$ be the usual involutive automorphism of~$B^n$ sending the origin in~$z_0/R$ (see [Ru, Section~2.2]), and let $\Phi_{z_0}\colon B^n\to B$ be given by $\Phi_{z_0}=R\phi_{z_0/R}$; in particular,
$\Phi_{z_0}$ is a biholomorphism with $\Phi_{z_0}(O)=z_0$, and thus
$\Phi_{z_0}\bigl(B_{B^n}(O,r)\bigr)=B_B(z_0,r)$. Furthermore (see [Ru, Theorem~2.2.6])
$$
|\hbox{Jac}_\R\Phi_{z_0}(z)|=R^{2n}\left({R^2-\|z_0\|^2\over
|R-\langle z,z_0\rangle|^2}\right)^{n+1}\ge {R^{n-1}\over 4^{n+1}}
\, d(z_0,\de B)^{n+1}\;,
$$
where $\hbox{Jac}_\R\Phi_{z_0}$ denotes the (real) Jacobian determinant of~$\Phi_{z_0}$.
It follows that
$$
\int_{B(z_0,r)}\chi\,d\nu=\int_{B_{B^n}(O,r)}(\chi\circ\Phi_{z_0})
|\hbox{\rm Jac}_\R\,\Phi_{z_0}|\,d\nu\ge {R^{n-1}\over 4^{n+1}}\,
d(z_0,\de B)^{n+1}\int_{B_{B^n}(O,r)}
(\chi\circ\Phi_{z_0})\,d\nu\;.
$$
Using [Ru, 1.4.3 and 1.4.7.(1)] we obtain
$$
\int_{B_{B^n}(O,r)}
(\chi\circ\Phi_{z_0})\,d\nu=2n
\int_{\de B^n}d\sigma(x){1\over 2\pi}
\int_0^r\int_0^{2\pi}\chi\circ\Phi_{z_0}(te^{i\theta}x)t^{2n-1}dt\,d\theta\;,
$$
where $\sigma$ is the area measure on~$\de B^n$ normalized so that $\sigma(\de B^n)=1$.
Now, $\zeta\mapsto \chi\circ\Phi_{z_0}(\zeta x)$ is subharmonic 
on~$r\Delta=\{|\zeta|<r\}\subset\C$ for any $x\in\de B^n$. Therefore [H\"o, Theorem~1.6.3] yields
$$
{1\over 2\pi}\int_0^r\int_0^{2\pi}\chi\circ\Phi_{z_0}(te^{i\theta}x)t^{2n-1}dt\,d\theta
\ge \chi(z_0)\int_0^r t^{2n-1}\,dt={1\over 2n}r^{2n}\chi(z_0)\;.
$$
So 
$$
\int_{B_{B^n}(O,r)}(\chi\circ\Phi_{z_0})\,d\nu\ge
r^{2n}\chi(z_0)\;,
$$
and the assertion follows.\qedn

%
%

Then:

\newthm Corollary \LSRdue: Let $D\subset\subset\C^n$ be a strongly pseudoconvex bounded  domain, and $r\in(0,1)$. Then there exists a $C_{3,r}>0$ depending on~$r$ such that
$$
\forevery{z_0\in D} 
\chi(z_0)\le {C_{3,r}\over \nu\bigl(B_D(z_0,r)\bigr)}\int_{B_D(z_0,r)}\chi\,d\nu
$$
for all non-negative plurisubharmonic functions $\chi\colon D\to\R^+$.

\pf
Since $D$ has smooth boundary, there exists a radius $\rho>0$ such that for every $x\in\de D$ the euclidean ball $B_x(\rho)$ internally tangent to $\de D$ at $x$ is contained in $D$. 

Let $z_0\in D$. If $d(z_0,\de D)\le\rho$, let $x\in\de D$ be such that
$d(z_0,\de D)=\|z_0-x\|$; in particular, $z_0$ belongs to the
ball $B=B_x(\rho)\subset D$. If $d(z_0,\de D)>\rho$, let $B\subset D$ be the 
Euclidean ball of center~$z_0$ and radius $d(z_0,\de D)$. In both cases
we have $d(z_0,\de D)=d(z_0,\de B)$; moreover,
the decreasing property of the Kobayashi distance yields $B_D(z_0,r) \supseteq 
B_B(z_0,r)$ for all $r\in(0,1)$.

Let $\chi$ be a non-negative plurisubharmonic function. Then
Lemmas~\ellisse\ and~\sei\ imply
$$
\eqalign{
\int_{B_D(z_0,r)}\chi\, d\nu&\ge \int_{B_B(z_0,r)}\chi\, d\nu \ge {\rho^{n-1}
\over 4^{n+1}}\,{r^{2n}d(z_0,\de D)^{n+1}\over
\nu\bigl(B_D(z_0,r)\bigr)}\,\nu\bigl(B_D(z_0,r)\bigr)
 \chi(z_0)\cr
&\ge {\rho^{n-1}
\over 4^{n+1}C_{1,r}}\,\nu\bigl(B_D(z_0,r)\bigr)\chi(z_0)\;,\cr}
$$
and we are done.
\qedn

\newrem If $D$ is strongly convex then one can take $C_{3,r}=C_3(1-r^2)^{-(n+1)}$, where
$C_3>0$ is independent of~$r$.

In a similar way we get another useful estimate:

\newthm Corollary \due: Let $D\subset\subset\C^n$ be a strongly pseudoconvex bounded  domain. Given $r\in(0,1)$, set $R={1\over 2}(1+r)\in(0,1)$. Then there exists a $K_r>0$ depending on~$r$ such that
$$
\forevery{z_0\in D\;\forall z\in B_D(z_0,r)} \chi(z)\le 
{K_r\over\nu\bigl(B_D(z_0,r)\bigr)}\int_{B_D(z_0,R)}\chi\,d\nu
$$
for every nonnegative plurisubharmonic function $\chi\colon D\to\R^+$.

\pf Let $r_1={1\over 2}(1-r)$; the triangle inequality implies that
$z\in B_D(z_0,r)$ yields $B_D(z,r_1)\subset B_D(z_0,R)$. Corollary~\LSRdue\
then implies
$$
\eqalign{
\chi(z)&\le {C_{3,r_1}\over\nu(B_D(z,r_1))}\int_{B_D(z,r_1)}\chi\,d\nu
\le {C_{3,r_1}\over\nu(B_D(z,r_1))}\int_{B_D(z_0,R)}\chi\,d\nu\cr
& = C_{3,r_1}{\nu(B_D(z_0,r))\over \nu(B_D(z,r_1))}
\cdot{1\over \nu(B_D(z_0,r))}\int_{B_D(z_0,R)}\chi\,d\nu
\cr}
$$
for all $z\in B_D(z_0,r)$.
Now Lemmas~\cinque\ and~\sei\ yield
$$
{\nu(B_D(z_0,r))\over \nu(B_D(z,r_1))}\ \leq\ {C_{1,r} C_2^{n+1}\over c_1(1-r)^{n+1}r_1^{2n}}
$$
and so
$$
\chi(z)\le {C_{3,r_1}C_{1,r}C_2^{n+1}\over c_1r_1^{2n}(1-r)^{n+1}}\, {1 \over \nu\bigl(B_D(z_0,r)
\bigr)}\int_{B_D(z_0,R)}\chi\,d\nu\;.
$$
\qedn


\smallsect 2. Carleson measures

Let $D\subset\subset\C^n$ be a strongly pseudoconvex bounded  domain in
$\C^n$. Given $0< p<+\infty$ the {\sl Bergman space}~$A^p(D)$ of~$D$ is the Banach space of holomorphic $L^p$-functions on~$D$, that is 
$A^p(D)=L^p(D)\cap\ca O(D)$, endowed with the $L^p$-norm 
$$
\|f\|_p^p=\int_D |f(z)|^p \,d\nu\;,
$$
where $\nu$ is the Lebesgue measure normalized so that $\nu(B^n)=1$. 

A finite positive Borel measure $\mu$ on~$D$ is said to be a {\sl Carleson measure}
of~$A^p(D)$ if there exists $C_p>0$ such that
$$
\forevery{f\in A^p(D)} \int_D|f(z)|^p\,d\mu\le C_p\|f\|_p^p\;.
$$
As explained in the introduction, our aim is to give a characterization
of Carleson measures involving the Bergman kernel of $D$.

Let $K\colon D\times D\to\C$ be the {\sl Bergman kernel} of~$D$ (see, e.g., [Kr, Section~1.4]);
it has the reproducing property
$$
\forevery{f\in A^2(D)\;\forall z\in D}f(z)=\int_D K(z,\zeta)f(\zeta)\,d\nu\;.
$$
Since $K(\cdot,\zeta)=\bar{K(\zeta,\cdot)}\in A^2(D)$, in particular we have
$$
K(z,z)=\int_D |K(z,\zeta)|^2\,d\nu(\zeta)=\|K(z,\cdot)\|_2^2\;.
$$

For each $z_0\in D$ let $k_{z_0}\in A^2(D)$ be the {\sl normalized Bergman kernel} given by
$$
k_{z_0}(z)={K(z,z_0)\over\|K(\cdot,z_0)\|_2}=
{K(z,z_0)\over\sqrt{K(z_0,z_0)}}\;;
$$
clearly, $\|k_{z_0}\|_2=1$. The {\sl Berezin transform} $B\mu$ of a finite measure $\mu$ on $D$ is 
the function given by
$$
B\mu(z)=\int_D|k_z(\zeta)|^2\,d\mu(\zeta)
$$
for all $z\in D$. 

We recall the following estimate from above on the Bergman kernel:

\newthm Lemma \tre: Let $D\subset\subset\C^n$ be a strongly
pseudoconvex bounded  domain. Then there exists $C_4>0$ such that
$$
\forevery{z_0\in D} |K(z_0,z_0)|\le{C_4\over d(z_0,\de D)^{n+1}}\;.
$$

\pf It follows immediately from [H\"o, Theorem~3.5.1] or from [R, p.~186].
\qedn

\Nota According to Range, it holds for any smooth domain.

Our next result is an estimate from below on the Bergman kernel, valid close enough to
the boundary.

\newthm Lemma \quattro: Let $D\subset\subset\C^n$ be a strongly
pseudoconvex bounded  domain. Then
for every $r\in(0,1)$ there exist $c_{5,r}>0$ and $\delta_r>0$ such that 
if $z_0\in D$ satisfies $d(z_0,\de D)<\delta_r$ then
$$
\forevery{z\in B_D(z_0,r)} |K(z,z_0)|\ge {c_{5,r}\over 
d(z_0,\de D)^{n+1}}\;.
$$

\pf Let $\psi$ be a defining function for~$D$; in particular, there are $C_6$,~$c_6>0$
such that
$$
c_6\, d(z,\de D)\le |\psi(z)|\le C_6\, d(z,\de D)
\neweq\eqdpsiapp
$$
in a neighbourhood of~$\bar D$. 

The main Theorem~2 in [F] implies that there is $\eta>0$
so that
$$
|K(z,z_0)|\ge {c\over\bigl(\psi(z)+\psi(z_0)+\rho(z,z_0)\bigr)^{n+1}}
$$
for a suitable constant $c>0$ as soon as $d(z_0,\de D)+d(z,\de D)+
\|z-z_0\|<\eta$, where 
$$
\rho(z,z_0)=\|z-z_0\|^2+\bigl|\de\psi_{z_0}(z-z_0)|\;.
$$
Put
$$
\delta_r={\eta\over 3}\min\left\{1, {1-r\over C_2},{c_{2,r}\eta\over3}\right\}\;;
$$
then Lemmas~\cinque\ and \quattrobis\ imply that $d(z_0,\de D)<\delta_r$
yields $d(z_0,\de D)+d(z,\de D)+\|z-z_0\|<\eta$ for all~$z\in B_D(z_0,r)$.
Using Lemmas~\cinque\ and~\quattrobis\ and \eqdpsiapp\ we then get
that if $d(z_0,\de D)<\delta_r$ then
$$
\eqalign{
|K(z,z_0)|&\ge {c\over \psi(z_0)^{n+1}}\left(
{1\over 1+{\psi(z)\over\psi(z_0)}+{\rho(z,z_0)\over\psi(z_0)}}
\right)^{n+1}\cr
&\ge {c(1-r)^{n+1}\over C_6^{n+1}d(z_0,\de D)^{n+1}}
\left({c_6\over c_6+C_2C_6+1/c_{2,r}} \right)^{n+1}
\cr}
$$
for all $z\in B_D(z_0,r)$, and we are done.\qedn

\newrem When $D$ is strongly convex we can take $c_{5,r}=c_5(1-r^2)^{n+1}$ with
$c_5>0$ independent of~$r$.

As a corollary we get a crucial estimate from below for the normalized Bergman kernel:

\newthm Corollary \piu: Let $D\subset\subset\C^n$ be a strongly
pseudoconvex bounded domain. Then for every $r\in(0,1)$ there exist $c_{7,r}>0$ and 
$\delta_r>0$ such that 
if $z_0\in D$ satisfies $d(z_0,\de D)<\delta_r$ then
$$
\forevery{z\in B_D(z_0,r)}|k_{z_0}(z)|^2\ge {c_{7,r}\over d(z_0,\de D)^{n+1}}\;.
$$

\pf It follows from Lemmas~\tre\ and~\quattro, with $c_{7,r}=c_{5,r}^2/C_4$.\qedn

\newrem When $D$ is strongly convex we can take $c_{7,r}=c_7(1-r^2)^{2(n+1)}$ with
$c_7>0$ independent of~$r$.

%
%
%

Now we can finally prove the promised characterization for the Carleson measures of $A^p(D)$:

\newthm Theorem \Carleson: Let $\mu$ be a finite positive Borel measure on a strongly pseudoconvex bounded  domain $D\subset\subset\C^n$. Then the following statements are equivalent:
\smallskip
\itm{(i)} $\mu$ is a Carleson measure of $A^p(D)$ for some (and hence all) $p\in(0,+\infty)$;
\itm{(ii)} the Berezin transform of $\mu$ is a bounded function;
\itm{(iii)} for every $r\in(0,1)$ there exists $C_r>0$ such that 
$\mu\bigl(B_D(z_0,r)\bigr)\le C_r\nu\bigl(B_D(z_0,r)\bigr)$ for all $z_0\in D$;
\itm{(iv)} for some $r\in(0,1)$ there exists $C_r>0$ such that 
$\mu\bigl(B_D(z_0,r)\bigr)\le C_r\nu\bigl(B_D(z_0,r)\bigr)$ for all $z_0\in D$.

\pf (i) $\Longrightarrow$(ii): By [CM] we can assume that $\mu$ is a Carleson measure of~$A^2(D)$. Then
$$
B\mu(z_0)=\int_D|k_{z_0}(z)|^2\,d\mu(z)\le C\|k_{z_0}\|_2^2=C
$$
for a suitable $C>0$, and $B\mu$ is bounded.

(ii)$\Longrightarrow$(iii): Fix $r\in(0,1)$. If $d(z_0,\de D)\ge\delta_r$,
where $\delta_r>0$ is given by Lemma~\quattro, then using
Lemma~\sei\ we get
$$
\mu\bigl(B_D(z_0,r)\bigr)\le\mu(D)\le {\mu(D)\over c_1r^{2n}\delta_r^{n+1}}
\nu\bigl(B_D(z_0,r)\bigr)
$$
as desired.

Assume now $d(z_0,\de D)<\delta_r$. Since the Berezin transform is bounded,
there exists $C_8>0$ independent of~$z_0$ and~$r$ such that
$$
\int_{B_D(z_0,r)}|k_{z_0}(z)|^2\,d\mu(z)\le B\mu(z_0)\le C_8\;.
$$
Hence Corollary~\piu\ yields
$$
{c_{7,r}\over d(z_0,\de D)^{n+1}}\mu\bigl(
B_D(z_0,r)\bigr)\le C_8\;.
$$
Recalling Lemma~\sei\ we get
$$
\mu\bigl(B_D(z_0,r)\bigr)\le {C_8\over c_{7,r}}\,d(z_0,\de D)^{n+1}\le {C_8\over c_1c_{7,r} r^{2n}}
\,\nu\bigl(B_D(z_0,r)\bigr)
$$
and we are done in this case too.

(iii)$\Longrightarrow$(iv): obvious.

(iv)$\Longrightarrow$(i): It follows from Lemmas~\sei, \cinque, Corollary~\LSRdue\ and [Lu]; for the sake of completeness we give here
a slightly different proof. Let $\{z_k\}$ be the sequence
given by Lemma~\uno. Clearly
$$
\int_D |f(z)|^p\,d\mu(z)\le \sum_{k=1}^\infty\int_{B_D(z_k,r)}|f(z)|^p\,
d\mu(z)
$$
for all $f\in A^p(D)$. Since $|f|^p$ is plurisubharmonic and nonnegative, Corollary~\due\ and (v) yields
$$
\eqalign{
\int_{B_D(z_k,r)}|f(z)|^p\,d\mu(z)&\le {K_r\over\nu\bigl(B_D(z_k,r)\bigr)}\int_{B_D(z_k,r)}d\mu(z)
\int_{B_D(z_k,R)}|f(\zeta)|^p\,d\nu(\zeta)\cr
&=K_r{\mu\bigl(B_D(z_k,r)\bigr)\over\nu\bigl(B_D(z_k,r)\bigr)}
\int_{B_D(z_k,R)}|f(\zeta)|^p\,d\nu(\zeta)\cr
&\le K_rC_r\int_{B_D(z_k,R)}|f(\zeta)|^p\,d\nu(\zeta)\;,
\cr}
$$
where $R={1\over2}(1+r)$.
Hence
$$
\int_D |f(z)|^p\,d\mu(z)\le K_rC_r\sum_{k=1}^\infty\int_{B_D(z_k,R)}|f(\zeta)|^p\,d\nu(\zeta)\le K_rC_r m\|f\|_p^p\;,
\neweq\eqdfin
$$
where $m$ is given by Lemma~\uno, and so $\mu$ is a Carleson measure of~$A^p(D)$.\qedn

\newrem Notice that \eqdfin\ says that if $\mu$ is a Carleson measure then we can find
$C>0$ so that
$$
\int_D |f(z)|^p\,d\mu(z)\le C \|f\|_p^p
$$
for all $f\in A^p(D)$ and all $p\in(0,+\infty)$; in other words, the constant $C$ is independent of~$p$.

\smallsect 3. Uniformly discrete sequences

Let $(X,d)$ be a metric space.  A sequence $\Gamma=\{x_j\}\subset
X$ of points in $X$ is {\sl uniformly discrete} if there exists $\delta>0$
such that $d(x_j,x_k)\ge\delta$ for all $j\ne k$. In this case 
$\inf\limits_{j\ne k} d(x_j,x_k)$ is the
{\sl separation constant} of~$\Gamma$. Furthermore, given $x_0\in X$, $r>0$ and
a subset $\Gamma\subset X$, we shall denote by
$N(x_0,r,\Gamma)$ the number of points of~$\Gamma$ contained in the ball
of center~$x_0$ and radius $r$. 

We can use $N(x_0,r,\Gamma)$ to detect finite unions
of uniformly discrete sequences:

\newthm Lemma \otto: Let $X$ be a metric space, and $\Gamma=\{x_n\}_{n\in\N}\subset X$ a sequence in $X$. If there are $N\ge 1$ and $r>0$ such that $N(x,r,\Gamma)
\le N$ for all $x\in X$, then $\Gamma$ is the union of at most $N$ uniformly discrete sequences.

\pf We shall define $N$ disjoint uniformly discrete sequences $\Gamma_0,\ldots,\Gamma_{N-1}\subset\Gamma$ so that $\Gamma=\Gamma_0\cup\cdots\cup\Gamma_{N-1}$. To do so,
we start with $\Gamma_0=\cdots=\Gamma_{N-1}=\void$ and, arguing by induction on~$n$,
we shall put each~$x_n$ in a specified~$\Gamma_j$. As a matter of notation, if $x_i\in
\Gamma_j$ we shall write $m(x_i)=j$, and we shall denote by $B(x,r)$ the metric ball of
center~$x$ and radius~$r$.

Put $x_0\in\Gamma_0$. Assume we have already defined
$m(x_i)$ for $i\le n$, and consider $x_{n+1}$. By assumption, $\Gamma\cap B(x_{n+1},r)$ contains at most $N$ points, one of which is $x_{n+1}$. Hence $\{x_0,\ldots,x_n\}\cap B(x_{n+1},r)$ contains at most $N-1$ points, and we can define
$$
m(x_{n+1})=\min\{i\in\{0,\ldots,N-1\}\ |\ i\neq m(x_j) \hbox{ for all $0\le j\leq n$ such that $x_j\in B(x_{n+1},r)$}\}\;.
$$
In this way $d(x_{n+1},x_j)\ge r$ for all $x_j\in\Gamma_{m(x_{n+1})}$ with $j<n+1$. 

It now follows easily that $\Gamma_0,\ldots,\Gamma_{N-1}$ are uniformly discrete sequences with separation constant at least $\delta=r$, because by construction if $x_h$,~$x_k\in\Gamma_j$ with $h>k$ we have $d(x_h,x_k)\ge r$. \qedn



We are now able to prove that uniformly discrete sequences give examples of
Carleson measures:

\newthm Theorem \dtre: Let $D\subset\subset\C^n$ be a strongly pseudoconvex bounded  domain, considered as a metric space with the distance~$\rho_D=\tanh k_D$. Let $\Gamma=\{z_j\}_{j\in\N}$ be a sequence
in~$D$. Then the following statements are equivalent:
\smallskip
\itm{(i)} $\Gamma$ is a finite union of uniformly discrete sequences;
\itm{(ii)} $\sup\limits_{z_0\in D}N(z_0,r,\Gamma)<+\infty$ for some $r\in(0,1)$;
\itm{(iii)} $\sup\limits_{z_0\in D}N(z_0,r,\Gamma)<+\infty$ for all $r\in(0,1)$;
\itm{(iv)}  there exists $p\in(0,+\infty)$ such that $\sum\limits_{z_j\in\Gamma} d(z_j,\de D)^{n+1}\delta_{z_j}$ is a Carleson measure of $A^p(D)$,
where $\delta_{z_j}$ is the Dirac measure in~$z_j$;
\itm{(v)} there exists $C_9>0$ such that for all $p\in(0,+\infty)$ we have
$$
\forevery{f\in A^p(D)} \sum_{z_j\in\Gamma} d(z_j,\de D)^{n+1}|f(z_j)|^p
\le C_9\|f\|_p^p\;;
$$
\indent in particular, $\sum\limits_{z_j\in\Gamma} d(z_j,\de D)^{n+1}\delta_{z_j}$ is a Carleson measure of $A^p(D)$
for all $p\in(0,+\infty)$.

\pf (v)$\Longrightarrow$(iv): Obvious.

(iv)$\Longrightarrow$(iii): Fix $r\in(0,1)$, and let $\delta_r>0$ be given
by Lemma~\quattro. By Lemma~\cinque, if $d(z_0,\de D)\ge\delta_r$ then $z\in B_D(z_0,r)$ implies
$d(z,\de D)\ge {1\over C_2}(1-r)\delta_r$. Using (iv) it is easy to see
that only a finite number of~$z_j\in\Gamma$ can have $d(z_j,\de D)\ge
{1\over C_2}(1-r)\delta_r$; therefore to get (iii) it suffice to prove
that the supremum is finite when $d(z_0,\de D)<\delta_r$.

Given $z_0\in D$ with $d(z_0,\de D)<\delta_r$, Corollary~\piu\ and Lemma~\cinque\ yield
$$
\forevery{z\in B_D(z_0,r)}
d(z,\de D)^{n+1}|k_{z_0}(z)|^2\ge   {c_{7,r}\over C_2^{n+1}}\,(1-r)^{n+1}\;.
$$
By [CM] we can assume $p=2$; hence
$$
N(z_0,r,\Gamma) \le{C_2^{n+1}\over c_{7,r}(1-r)^{n+1}}
\sum_{z\in B_D(z_0,r)\cap\Gamma}d(z,\de D)^{n+1}|k_{z_0}(z)|^2
\le {C_2^{n+1}C\over c_{7,r}(1-r)^{n+1}}\|k_{z_0}\|_2^2=
{C_2^{n+1}C\over c_{7,r}(1-r)^{n+1}}
$$
for a suitable $C>0$, as desired.

(iii)$\Longrightarrow$(ii): Obvious.

(ii)$\Longrightarrow$(i): Lemma~\otto.

(i)$\Longrightarrow$(v): Clearly it suffices to prove the assertion when $\Gamma$ is a single
uniformly discrete sequence. Let $\delta>0$ be the separation constant
of~$\Gamma$, and put $r=\delta/2$.
By the triangle inequality, the Kobayashi balls $B_D(z_j,r)$ are
pairwise disjoint. Hence
$$
\int_D|f(z)|^p\,d\nu\ge\sum_{z_j\in\Gamma}\int_{B_D(z_j,r)}|f(z)|^p
\,d\nu\;.
$$
Now, $|f|^p$ is plurisubharmonic and nonnegative; hence Corollary~\LSRdue\ 
and Lemma~\sei\ yield
$$
\int_{B_D(z_j,r)}|f(z)|^p\,d\nu\ge {1\over C_{3,r}}\nu(B_D(z_j,r))|f(z_j)|^p\ge {c_1\over C_{3,r}}r^{2n}d(z_j,\de D)^{n+1}|f(z_j)|^p\;,
$$
for all $z_j\in\Gamma$.
Setting $C_9={C_{3,r}\over c_1 r^{2n}}$, the assertion
follows.\qedn

Now we would like to study how fast a uniformly discrete sequence can escape to the
boundary. A first result in this vein is an immediate corollary of the previous theorem:

\newthm Corollary \duno: Let $\Gamma=\{z_j\}\subset D$ be a uniformly discrete
sequence in a strongly pseudoconvex bounded domain $D\subset\subset
\C^n$. Then
$$
\sum_{z_j\in\Gamma} d(z_j,\de D)^{n+1}<+\infty\;.
$$

\pf It suffices to take $f\equiv 1$ in Theorem~\dtre.(v).\qedn



The next result gives a worse estimate, but valid in any
hyperbolic domain with finite Euclidean volume. 

\newthm Proposition \dqua: Let $D\subset\C^n$ be a hyperbolic domain with finite Euclidean
volume, endowed with the distance $\rho_D=\tanh k_D$.
Let $\Gamma=\{z_j\}\subset D$ be an uniformly discrete sequence with $d(z_j,\de D)<1$ for all~$z_j\in\Gamma$. Then
$$
\sum_{z_j\in\Gamma} d(z_j,\de D)^{2n}\, h\left(-{1\over
\log d(z_j,\de D)}\right)<+\infty
$$
for any increasing function $h\colon\R^+\to\R^+$ such that
$$
\sum_{m=1}^{+\infty}h\left({1\over m}\right)<+\infty\;.
$$

\pf For $m\in\N$ put
$$
\Omega_m=\left\{z\in D\biggm| m<\log{1\over d(z,\de D)}\le m+1\right\}=
\{z\in D\mid e^{-(m+1)}\le d(z,\de D)< e^{-m}\}\;.
$$
Set $\Gamma_m=\Gamma\cap\Omega_m$; by assumption, $\Gamma=\bigcup\limits_{m=0}^{+\infty}\Gamma_m$. Let $\delta>0$ be the separation constant
of~$\Gamma$, and put $r=\delta/2$. The Kobayashi balls $B_D(z_j,r)$
centered in points of~$\Gamma$ are pairwise disjoint; in particular,
$$
\sum_{z_j\in\Gamma_m}\nu\bigl(B_D(z_j,r)\bigr)\le \nu(D)\;.
\neweq\eqdpiu
$$
%
If $z\in D$ then the euclidean
ball $B$ of center~$z$ and radius $d(z,\de D)$ is contained in~$D$. But then
$B_D(z,r)$ contains the Kobayashi ball~$B_B(z,r)$ of $B$
centered in~$z$ and of radius~$r$, which is an euclidean
ball of radius $r d(z,\de D)$. Thus 
$$
\nu\bigl(B_D(z_j,r)\bigr)\ge r^{2n}d(z_j,\de D)^{2n}\;.
$$
Let $N_m$ be the cardinality of~$\Gamma_m$. Then \eqdpiu\ yields
$$
N_m \le \nu(D)r^{-2n}e^{2n(m+1)}=C_9 e^{2nm}\;,
$$
for a suitable constant $C_9>0$ independent of~$m$.
%
Therefore
$$
\eqalign{
\sum_{z_j\in\Gamma} d(z_j,\de D)^{2n}\, h\left(-{1\over
\log d(z_j,\de D)}\right)&=\sum_{m=0}^{+\infty}
\sum_{z_j\in\Gamma_m}d(z_j,\de D)^{2n}\, h\left(-{1\over
\log d(z_j,\de D)}\right)\cr
&\le \sum_{z_j\in\Gamma_0}d(z_j,\de D)^{2n}\, h\left(-{1\over
\log d(z_j,\de D)}\right)+
\sum_{m=1}^{+\infty} N_m e^{-2mn}h\left({1\over m}\right)\cr
&\le C_{10}+C_9\sum_{m=1}^{+\infty} h\left({1\over m}\right)<+\infty\;.
\cr}
$$
\qedn

\newrem If $D\subset\subset\C^n$ is bounded and $\Gamma=\{z_j\}\subset D$ 
is an uniformly discrete (with respect to any distance inducing the natural topology of $D$) sequence, then $d(z_j,\de D)<1$ for all but a
finite number of elements of~$\Gamma$.

To get the sharp estimate valid in strongly pseudoconvex bounded  domains,
we replace the Euclidean measure by the 
Eisenman-Kobayashi invariant measure (but see Remark~3.2 below).

For simplicity, let us recall the definition of the Eisenman-Kobayashi measure in domains of~$\C^n$ only. Let $D\subset\C^n$ be a domain;
the {\sl Eisenman-Kobayashi density}~$K_D\colon D\to\R^+$ is
given by
$$
K_D(z_0)=\inf\bigl\{|\hbox{Jac}_\R f(O)|^{-1}\bigm| 
f\colon B^n\to D\hbox{ holomorphic, } f(O)=z_0, \hbox{ $df_O$ invertible}
\bigr\}\;.
$$
It is not difficult to prove that $K_D$ is upper 
semicontinuous (see, e.g., [A, Proposition~2.3.37]); then
the {\sl Eisenman-Kobayashi measure} of~$D$ is just $\tilde\kappa_D=K_D\,\nu$.
When $D=B^n$, the Eisenman-Kobayashi measure coincides with the
Bergman (or hyperbolic) volume (see, e.g., [A, Proposition~2.3.36]). 

The next lemma contains an estimate on the Eisenman-Kobayashi measure
of Kobayashi balls.

\newthm Lemma \nove: Let $D\subset\subset\C^n$ be a 
strongly pseudoconvex bounded  domain,. Then there exist $0<c_{11}$, $C_{11}$ such that
$$
\forevery{z_0\in D\;\forall r\in(0,1)}\qquad\qquad\quad c_{11}r^{2n}(1-r)^{n+1}\le
\tilde\kappa_D\bigl(B_D(z_0,r)\bigr)\le {C_{11}\over d(z_0,\de D)^n(1-r)^n}\;.
$$  

\pf The estimates in [M] show that there are two constants
$c_{12}$, $C_{12}>0$ such that
$$
\forevery{z\in D}
{c_{12}\over d(z,\de D)^{n+1}}\le K_D(z)\le 
{C_{12}\over d(z,\de D)^{n+1}}\;.
$$
In particular, Lemmas~\cinque\ and~\sei\ yield
$$
\eqalign{
\tilde\kappa_D\bigl(B_D(z_0,r)\bigr)&=\int_{B_D(z_0,r)}K_D(z)\,d\nu\ge
c_{12}\int_{B_D(z_0,r)}{1\over d(z,\de D)^{n+1}}\,d\nu\ge
{c_{12}(1-r)^{n+1}\over C_2^{n+1}d(z_0,\de D)^{n+1}}\nu\bigl(B_D(z_0,r)\bigr)\cr
&\ge  {c_{12} c_1\over C_2^{n+1}} r^{2n}(1-r)^{n+1}\;.
\cr}
$$
For the upper estimate, notice that since $\de D$ is  and bounded there is an $\eps>0$ such that, setting 
$$
U_\eps = D\setminus K_\eps = \left\{z\in D\, |\, d(z,\de D)<\eps\right\}\;,
$$
then there exists a smooth diffeomorphism $\Psi\colon\overline U_\eps\to [0,\eps]\times\de D$ with $\Psi^{-1}(\{t\}\times\de D)=\{z\in D\mid d(z,\de D)=t\}$ for all $t\in(0,\eps]$. 

Now Lemma~\cinque\ implies
$$
B_D(z_0,r)\subset \left\{z\in D\, |\, d(z,\de D)\geq l(r)\right\}\;,
$$
where $l(r)={1\over C_2}(1-r)d(z_0,\de D)$. Thus, denoting by $\sigma$
the usual $(2n-1)$-measure on~$\de D$, and by $\lambda$
the product measure on $[0,\eps]\times\de D$, and using the 
compactness of $\overline U_\eps$, we have
$$
\eqalign{
\tilde\kappa_D\bigl(B_D(z_0,r)\bigr)&=\int_{B_D(z_0,r)}K_D(z)\,d\nu
\le C_{12}\int_{B_D(z_0,r)}{1\over d(z,\de D)^{n+1}}\,d\nu\cr
&\le C_{12}\int_{K_\epsilon}{1\over d(z,\de D)^{n+1}}\,d\nu +C_{12} \int_{B_D(z_0,r)\cap U_\epsilon}{1\over d(z,\de D)^{n+1}}\,d\nu\cr
&\le C_{13} + C_{14} \int_{\Psi(B_D(z_0,r)\cap U_\epsilon)}{1\over t^{n+1}}\,d\lambda\le C_{13} + C_{14} \int_{[l(r),\eps]\times\de D}{1\over t^{n+1}}\,d\lambda\cr
&\le C_{13} + C_{14}\,\sigma(\de D) \left({1\over l(r)^n}-{1\over \epsilon^n}\right) \le {C_{15}\over l(r)^n }={C_2^n C_{15}\over d(z_0,\de D)^n (1-r)^n}\;,\cr}
$$
for suitable constants $C_{13}$, $C_{14}$, $C_{15}>0$ independent of~$z_0$
and~$r$.\qedn

\newrem We stated this lemma in terms of the Eisenman-Kobayashi measure
just to keep with the invariant approach of this paper; but for the purpose of the next
theorem any measure providing the same estimates would work. For instance,
we might use the measure $\mu=d(\cdot,\de D)^{-(n+1)}\nu$.

We are now able to prove the promised sharp estimate for strongly pseudoconvex domains:

\newthm Theorem \dqua: Let $D\subset\subset\C^n$ be a strongly pseudoconvex bounded domain, endowed with the distance $\rho_D=\tanh k_D$.
Let $\Gamma=\{z_j\}\subset D$ be an uniformly discrete sequence with $d(z_j,\de D)<1$ for all~$z_j\in\Gamma$. Then
$$
\sum_{z_j\in\Gamma} d(z_j,\de D)^{n}\, h\left(-{1\over\log
d(z_j,\de D)}\right)<+\infty
$$
for any increasing function $h\colon\R^+\to\R^+$ such that
$$
\sum_{m=1}^{+\infty}h\left({1\over m}\right)<+\infty\;.
$$

\pf Fix $z_0\in D$, and for $m\in\N$ set
$$
D_m=\left\{z\in D\biggm|{m\over 2}\le k_D(z_0,z)<{m+1\over 2}\right\}\;, 
$$
and let $\Gamma_m=D_m\cap \Gamma$. Let $\delta>0$ be the separation
constant of~$\Gamma$ (with respect to~$k_D$), and let $r=\tanh(\delta/2)$. 
Clearly, if $z_j\in\Gamma_m$ and $z\in B_D(z_j,r)$, then $k_D(z_0,z)<{1\over 2}(m+1+\delta)$. Since, as usual, the Kobayashi balls $B_D(z_j,r)$
are pairwise disjoint, using Lemma~\nove\ we get
$$
c_{11}r^{2n}(1-r)^{n+1} N_m\le\sum_{z_j\in\Gamma_m}\tilde\kappa_D\bigl(B_D(z_j,r)\bigr)
\le \tilde\kappa_D\bigl(B_D(z_0,R_m)\bigr)\le{C_{11}\over
d(z_0,\de D)^n(1-R_m)^n}\;,
$$
where $R_m=\tanh{m+1+\delta\over 2}$; hence there is $C_{16}>0$ (depending
on~$z_0$ and $\delta$ but not on~$m$) such that
$$
N_m\le C_{16} e^{mn}\;.
$$

Now, the estimates \equbb\ on the boundary behavior of the Kobayashi distance yield $C_{17}$, $C_{18}>0$ such that
$$
d(z,\de D)\le C_{17}\exp\bigl(-2 k_D(z_0,z)\bigr)
$$
and
$$
-{1\over\log d(z,\de D)}\le {1\over 2k_D(z_0,z)-C_{18}}
$$
as soon as $2k_D(z_0,z)>C_{18}$.
In particular, if $z_j\in\Gamma_m$ we have
$$
d(z_j,\de D)\le C_{17}e^{-m}\qquad\hbox{and}\qquad
-{1\over\log d(z,\de D)}\le{1\over m-C_{18}}\;.
$$
as soon as $m>C_{18}$.
Then
$$
\eqalign{
\sum_{z_j\in\Gamma} d(z_j,\de D)^{n}\,& h\left(-{1\over\log
d(z_j,\de D)}\right)=\sum_{m=0}^{+\infty}\sum_{z_j\in\Gamma_m}
d(z_j,\de D)^{n}\,h\left(-{1\over\log
d(z_j,\de D)}\right)\cr
&\le C_{19}+C_{17}^n\sum_{m=m_0}^{+\infty}N_m e^{-nm}h\left({1\over m-C_{18}}\right)
\le C_{19}+C_{17}^nC_{16}\sum_{m=0}^{+\infty}h\left({1\over m}\right)<+\infty\;,
\cr}
$$
where $m_0$ is the smallest integer greater than $C_{18}$.\qedn

In particular, taking $h(x)=e^{-1/x}$ we recover Corollary~\duno.

\newrem The statement of Theorem~\dqua\ is sharp. In fact, MacCluer [Mc] has
constructed for any increasing function $h\colon\R^+\to\R^+$ with
$\lim_{t\to 0^+}h(t)=0$ and $\sum_{m=1}^\infty h(1/m)=+\infty$
an uniformly discrete sequence $\Gamma=\{z_j\}$ in~$B^n$ such that
$$
\sum_{z_j\in\Gamma}d(z_j,\de B^n)^n h\left(-{1\over\log d(z_j,\de B^n)}
\right)=+\infty\;.
$$

\setref{DSV}

\beginsection References

\book A M. Abate: Iteration theory of holomorphic maps on taut manifolds! Mediterranean Press, Cosenza, 1989; see {\tt http://www.dm.unipi.it/\lower3pt\hbox{\tt \char"7E}abate/libri/libriric/libriric.html}

\art C L. Carleson: Interpolations by bounded analytic functions and the corona problem!
Ann. of Math.! 76 1962 547-559

\art CM J.A. Cima, P.R. Mercer: Composition operators between Bergman spaces 
on convex domains in $\C^n$! J. Operator Theory! 33 1995 363-369

\art CW J.A. Cima, W.R. Wogen: A Carleson measure theorem for the Bergman space on the ball!
J. Operator Theory! 7 1982 157-165

\art D P.L. Duren: Extension of a theorem of Carleson! Bull. Amer. Math. Soc.! 75 1969 143-146

\book DS P.L. Duren, A. Schuster: Bergman spaces! American Mathematical Society, Providence, RI, 2004

\coll DSV P.L. Duren, A. Schuster, D. Vukoti\'c: On uniformly discrete sequences in the disk!
Quadrature domains and applications! P. Ebenfeld, B. Gustafsson, D. Khavinson and M. Putinar eds., Birkh\"auser, Basel, 2005, pp. 131--150

\art DW P.L. Duren, R. Weir: The pseudohyperbolic metric and Bergman spaces in the ball!
Trans. Amer. Math. Soc.! 359 2007 63-76

\art F C. Fefferman: The Bergman kernel and biholomorphic mappings of pseudoconvex domains! 
Invent. Math.! 26 1974 1-65

\art H W.W. Hastings: A Carleson measure theorem for Bergman spaces! Proc. Amer. Math. Soc.! 52 1975 237-241


\book H\"o L. H\"ormander: An introduction to complex analysis in several variables!
North Holland, Amsterdam, 1973

\art JMT M. Jevti\'c, X. Massaneda, P.J. Thomas: Interpolating sequences for weighted Bergman 
spaces of the ball! Michigan Math. J.! 43 1996 495-517

\book JP M. Jarnicki, P. Pflug: Invariant distances and metrics in complex analysis! Walter de Gruyter \& co., Berlin, 1993

\book K S. Kobayashi: Hyperbolic complex spaces! Springer-Verlag, Berlin, 1998

\book Kr S. Krantz: Function theory of several complex variables! Wiley, New York, 1981

\art KM S.G. Krantz, D. Ma: Bloch functions on strongly pseudoconvex domains! Indiana 
Univ. Math. J.! 37 1988 145-163


\art Li H. Li: BMO, VMO and Hankel operators on the Bergman space of strictly pseudoconvex
domains! J. Funct. Anal.! 106 1992 375-408

\art Lu D. Luecking: A technique for characterizing Carleson measures on Bergman spaces! Proc. Amer. Math. Soc.! 87 1983 656-660

\art M D. Ma: Boundary behavior of invariant metrics and volume forms on strongly pseudoconvex domains! Duke Math. J.! 63 1991 673-697

\art Mc B.D. MacCluer: Uniformly discrete sequences in the ball! J. Math. Anal. Appl.! 318 2006 37-42 

\art Ma X. Massaneda: $A^{-p}$ interpolation in the unit ball! J. London Math. Soc.! 52 1995
391-401

\art O V.L. Oleinik: Embeddings theorems for weighted classes of harmonic and analytic functions!
J. Soviet Math.! 9 1978 228-243

\art OP V.L. Oleinik, B.S. Pavlov: Embedding theorems for weighted classes of harmonic and analytic functions! J. Soviet Math.! 2 1974 135-142

\book R R.M. Range: Holomorphic functions and integral representations in several complex variables! Springer-Verlag, Berlin, 1986

\book Ru W. Rudin: Function theory in the unit ball of $\C^n$! Springer-Verlag, Berlin, 1980

\bye